\documentclass[12pt]{iopart}
\usepackage{iopams}  
\begin{document}
% Journal identifier can be put here if required, e.g.
\jl{01}

\title{Random projections and the optimization
of an algorithm for phase retrieval}

\author{Veit Elser}

\address{Department of Physics,
Cornell University, Ithaca, NY 14853-2501, USA}

\begin{abstract}
Iterative phase retrieval algorithms typically employ projections onto constraint subspaces
to recover the unknown phases in the Fourier transform of an image, or, in the case of
x-ray crystallography, the electron density of a molecule. For a general class of algorithms,
where the basic iteration is specified by the {\em difference map},
solutions are associated with fixed points of the map, the attractive character of which
determines the effectiveness of the algorithm. The behavior of the difference map near fixed points
is controlled by the relative orientation of the tangent spaces of the two constraint subspaces
employed by the map. Since the dimensionalities involved are always large in practical
applications, it is appropriate to use random matrix theory ideas to analyze the average-case
convergence at fixed points. Optimal values of the $\gamma$ parameters of the
difference map are found which differ somewhat from the values previously obtained
on the assumption of orthogonal tangent spaces.
\end{abstract}


\submitted

\maketitle

\section{Introduction}
In two dimensional interferometric imaging, as well as diffraction from three dimensional
crystals and non-periodic objects, the Fourier modulus of an object is sampled
on a discrete and finite grid. To recover the object, additional {\it a priori} constraints
(object support, positivity, atomicity, etc.) are imposed in order to ``retrieve"
the corresponding phases, without which the Fourier synthesis of the object is
impossible. A natural setting for computations is the $N$ dimensional complex
Euclidean space $E^N$, where $N$ is the number of pixels (or voxels)
in the finite Fourier transform grid.
The most widely used phase retrieval algorithms are formulated
as maps $E^N\to E^N$, which, when used iteratively, impose the Fourier modulus data
as well as the {\it a priori} constraints in the course of recovering the object.

The work reported here concerns the convergence of the phase retrieval
algorithm based on the {\em difference map}\cite{elser}:
\begin{equation}\label{diffMapDefined}
\rho\mapsto D(\rho) = \rho + \beta \left( \Pi_1 \circ f_2 - \Pi_2 \circ f_1 \right)(\rho)\; .
\end{equation}
The action of $D$ is to add to the current object,
$\rho\in E^N$, a difference of projections, $\Pi_1$ and $\Pi_2$. In general, the action of a projection
on an object $\rho$ is the minimal modification of $\rho$ that restores a particular constraint, where
minimality is specified by the standard Euclidean norm in $E^N$. 
The two specific constraints considered below
are the Fourier modulus constraint, implemented by $\Pi_{\rm F}$, and a ``fixed" or ``atomic" support
constraint, implemented by $\Pi_{\rm S}$. The nonzero real parameter $\beta$ serves as the step
size of the iterations. Convergence of the difference map is crucially dependent on the
maps $f_1$ and $f_2$ with which the basic projections in (1) are composed:
\begin{equation}
f_i(\rho) = (1 + \gamma_i)\Pi_i(\rho) - \gamma_i\, \rho\qquad (i = 1, 2)\;.
\end{equation}
The structure of the maps $f_i$ is geometrical, taking $\rho$ to a general point parametrized by
$\gamma_i$ on the line determined by
$\rho$ and $\Pi_i(\rho)$. A previous analysis\cite{elser} of the convergence of the difference map,
greatly simplified by the assumption of orthogonality of the two constraint subspaces near
the solution, found
\numparts
\begin{eqnarray}
\gamma_1^{\rm opt} &= -1/\beta \\
\gamma_2^{\rm opt} &= 1/\beta\;,
\end{eqnarray}
\endnumparts
as optimal parameter values. The main result reported here is a modification of these values
for the more realistic case of constraint subspaces that are not perfectly orthogonal, locally,
in the vicinity of the solution. 

The presentation is organized as follows. Section 2 reviews the relationship between the fixed
points of the difference map and the solution to the associated phase retrieval problem. Explicit
expressions for the two chief projections and their linearization are derived in Section 3.
Convergence of the difference map for a general pair of linearized projections is studied in Section 4.
Traces encountered in Section 4 are averaged in a standard random matrix ensemble in Section 5, and
compared with averages over object ensembles appropriate for phase retrieval. 
Results are summarized in Section 6. 

\section{The difference map}
A special case of the difference map, the {\em hybrid input-output map}\cite{fienup1}, 
first appeared in the context of two dimensional phase retrieval with a fixed support constraint. 
The generality of the construction for an
arbitrary pair of projections\cite{elser, bauschke}, 
and the flexibility in choosing the parameters $\gamma_i$\cite{elser}, was only
noticed recently. With $\Pi_1 = \Pi_{\rm S}$ (support), $\Pi_2 = \Pi_{\rm F}$ (Fourier modulus),
the hybrid input-output map is obtained for the parameter values $\gamma_1 = -1$, $\gamma_2 = 1/\beta$
\cite{elser}.

A phase retrieval solution $\rho^{\rm sol}$ exists if and only if the difference map has
a fixed point $\rho^*$. This is an immediate consequence of the definition of the difference map. 
At a fixed point
$\rho^* = D(\rho^*)$, the difference of projections in (\ref{diffMapDefined}) vanishes, and we have
\begin{equation}\label{fixPtProj}
\Pi_1\circ f_2(\rho^*) = \Pi_2\circ f_1(\rho^*) := \rho^{\rm sol}\;.
\end{equation}
From this follows
\begin{equation}
\Pi_1(\rho^{\rm sol}) = \Pi_2(\rho^{\rm sol}) = \rho^{\rm sol}\;,
\end{equation}
showing that the object $\rho^{\rm sol}$ satisfies both constraints. Conversely, given $\rho^{\rm sol}$,
the set of fixed points is given by
\begin{equation}\label{fixedPts}
(\Pi_1\circ f_2)^{-1}(\rho^{\rm sol}) \cap (\Pi_2\circ f_1)^{-1}(\rho^{\rm sol})\;,
\end{equation}
a nonempty set since it trivially contains $\rho^{\rm sol}$.

To recover an object using the difference map one begins with an initial, arbitrary object $\rho(0)$, and generates iterates $\rho(n) = D^n(\rho(0))$ until the norm of the difference
\begin{equation}
\|\rho(n+1) - \rho(n)\| = \|\beta\left(\Pi_1 \circ f_2 - \Pi_2 \circ f_1\right)(\rho(n))\|
\end{equation}
becomes suitably small. Then, from the approximate fixed point $\rho^*\approx\rho(n)$, the solution 
$\rho^{\rm sol}$ is obtained using (\ref{fixPtProj}). There is no guarantee, that given an arbitrary
starting point $\rho(0)$, the iterates will converge on a fixed point. In particular, it may happen that
the linearization of the difference map in the neighborhood of a fixed point has unstable
directions. Avoiding this, and instead making fixed points maximally stable in an average sense,
is the primary focus of the present work.

The uniqueness of solutions is linked to the tightness of the constraints. In many applications the
two constraint subspaces $C_1$ and $C_2$ are unions of countably many submanifolds of $E^N$ and
their stable intersection properties ({\it i.e.} solutions) are functions of their dimensionalities. Thus if
$\dim(C_1) + \dim(C_2) > N$, then the set of solutions, $C_1\cap C_2$, has positive dimension and
there is no uniqueness. For phase retrieval to be well posed, the number of {\it a priori} constraints must 
be sufficiently large, and the dimensionality of the associated constraint subspace correspondingly small, so that only the empty intersection is stable. In this overconstrained situation a solution will reduce to
a single point (or low dimensional submanifold in the case of symmetries) since a solution is 
known to exist at the outset.

\section{Projections}
When the Fourier transform of an object is sampled on a regular grid (fundamentally unavoidable in crystallography,
an artifact of CCD array design otherwise), the object can be treated as living on a {\em periodic}
real-space grid. The fineness of the grid in real-space is determined by the range of the
Fourier-space sampling. A general point of the periodic real-space grid has the form
\begin{equation}
\bi{r} = (r_1, r_2, \ldots)\qquad r_i\in \mathbb{Z}_{m_i}/m_i\;,
\end{equation}
where the integers $m_i$ give the ranges of the Fourier-space sampling
for the different components of the Fourier-space basis, and $\mathbb{Z}_m$ denotes the set
of integers modulo $m$. Points in Fourier-space (relative to a basis which may be non-orthogonal)
have the form
\begin{equation}\label{FourierGrid}
 \bi{q} = (q_1, q_2, \ldots)\qquad q_i\in 2\pi\,\mathbb{Z}_{m_i}\;.
\end{equation}
In crystallography, $\bi{r}$ is called a ``fractional coordinate", while $\bi{q}$ is,
up to a factor of $2\pi$, a triplet of Miller indices. The periodicity of the Fourier-space grid
(\ref{FourierGrid}) is an artifact of the regular, grid-like sampling in real-space; this artificial
periodicity can be tolerated since the modulus of the Fourier transform of real objects
becomes negligible at the extremes of the range. With points in real- and Fourier-space defined
as above, the components of an arbitrary object in real-space, $\rho_{\bi{r}}$, and
Fourier-space, $\rho_{\bi{q}}$, are related by the unitary transformation (discrete Fourier
transform) with matrix elements
\begin{equation}
(\mathcal{F})_{\bi{q\, r}} = \frac{1}{\sqrt{N}}\exp\rmi \bi{q}\cdot\bi{r}\;,
\end{equation}
where $N = m_1\,m_2 \cdots$ is the total number of real- or Fourier-space samples.

\subsection{Support projection}
Object support constraints can take two forms, depending on the application. In
phase retrieval with non-periodic objects one normally has {\it a priori} knowledge of 
a fixed object support $S$, with the property $\rho_{\bi{r}} = 0$ if $\bi{r}\notin S$.
For example, relatively tight bounds on the set $S$ can be derived from the object's
autocorrelation\cite{fienup2}, as given directly by the Fourier transform of the squared Fourier modulus.
In crystallography, on the other hand, it is the {\em form} of the support that supplies the necessary
constraint. Given the number of atoms $M$, a valid support comprises a union of $M$ very
compact sets, say $3\times 3\times 3$ arrays of voxels. Although the locations of the
individual atomic supports is unknown, the small size of the combined support of all $M$
atoms is a very strong constraint.

In the case of a {\em fixed object support} the corresponding projection operator is linear and has the
following simple form:
\begin{equation}\label{suppProj}
(\Pi_{\rm S})_{\bi{r\, r^\prime}} = \sum_{\bi{s}\in S}\delta_{\bi{r\, s}}\delta_{\bi{s\, r^\prime}}\;.
\end{equation}
It is easily verified that given an arbitrary object $\rho$, the object $\Pi_{\rm S}(\rho)$ is the
nearest object (in the Euclidean sense) having support $S$. The corresponding constraint subspace
is a linear subspace of $E^N$ with dimensionality $|S|$.

The constraint subspace for an {\em atomic support} is, in contrast, a union of very many linear
subspaces, each one corresponding to a different arrangement of atoms. Thus one is really
working with a collection of supports $\mathcal{S}$, where the linear projection operator
(\ref{suppProj}) applies to each $S\in\mathcal{S}$. The actual projection operator is nonlinear,
returning the element of $\{\Pi_{\rm S}(\rho)\colon S\in\mathcal{S}\}$ which minimizes
$\|\Pi_{\rm S}(\rho) - \rho\|$. This complication, however, is not relevant to the local analysis
we perform below: only one element of $\mathcal{S}$ contains the solution (given a choice of
object origin and enantiomer) and the linear projection (\ref{suppProj}) can be used for just that
one.

\subsection{Fourier modulus projection}
We consider the situation where a large number of the object's Fourier moduli have
been measured to high precision. In this case, the moduli $F_{\bi{q}}$ in
\begin{equation}
\rho_{\bi{q}}^{\rm sol} = F_{\bi{q}}\exp{\rmi \phi_{\bi{q}}}\;,
\end{equation}
for a large subset $Q_{\rm data}$ of Fourier samples $\bi{q}$, can be treated as known quantities.
The corresponding projection operator is nonlinear but simple when expressed in Fourier-space:
\begin{equation}\label{FourierModProj}
\widetilde{\Pi}_{\rm F}(\rho_{\bi{q}}) = \cases{
F_{\bi{q}}\frac{\rho_{\bi{q}}}{|\rho_{\bi{q}}|} & for $\bi{q}\in Q_{\rm data}$ and $\rho_{\bi{q}}\neq 0$\\
F_{\bi{q}} & for $\bi{q}\in Q_{\rm data}$ and $\rho_{\bi{q}}= 0$\\
\rho_{\bi{q}} & otherwise.\\}
\end{equation}
On actual computing machinery, and random initialization, the arbitrary choice of phase made in the
second case above never comes up in practice. In crystallography, $Q_{\rm data}$ never includes
the origin and frequently also omits ``reflections" $\bi{q}$ close to the origin. Since in all applications
there is a systematic decay of moduli at large spatial frequencies, a significant fraction
of the samples near the corners of the Fourier grid (those farthest from the origin) 
will lack measured values. There is usually no harm in approximating
these small moduli by zero, since it is the
absolute error that matters in Fourier synthesis. 
Geometrically, Fourier modulus projection in Fourier-space corresponds to
projections of the complex numbers $\rho_{\bi{q}}$ onto circles with prescribed radii ($F_{\bi{q}}$),
for all $\bi{q}\in Q_{\rm data}$. A more refined approach, 
one which takes into account measurement errors and bounds on unmeasured moduli, would
involve projections, respectively, onto annuli and disks.
When considered together with support projection (\ref{suppProj}), it is necessary to use
the real-space form of Fourier modulus projection:
\begin{equation}\label{realSpaceFourierProj}
\Pi_{\rm F} = \mathcal{F}^{-1}\circ\widetilde{\Pi}_{\rm F}\circ\mathcal{F}\;.
\end {equation}

\subsection{Linearized projections}
To study convergence, the difference map and the projections from which it is built are
linearized about the solution $\rho^{\rm sol}$. In general, the linearization $\pi$ of a
projection $\Pi$ is defined by
\begin{equation}
\pi(\eta) :=\lim_{\epsilon\to 0}\frac{\Pi(\rho^{\rm sol}+\epsilon\eta) - \rho^{\rm sol}}{\epsilon}\;.
\end{equation}
Since object support projection (\ref{suppProj}) is already linear,
\begin{equation}\label{linearSuppProj}
\pi_{\rm S}(\eta) = \Pi_{\rm S}(\eta)\;.
\end{equation}
Linearized Fourier modulus projection (\ref{FourierModProj}) is diagonal in Fourier-space,
\begin{equation}
(\tilde{\pi}_{\rm F})_{\bi{q\, q^\prime}} = \delta_{\bi{q\, q^\prime}}\cdot
\cases{
{\textstyle \frac{1}{2}}-{\textstyle \frac{1}{2}}\exp{(\rmi 2\phi_{\bi{q}})}\; C&
for $\bi{q}\in Q_{\rm data}$ and $F_{\bi{q}}\neq 0$\\
0&
for $\bi{q}\in Q_{\rm data}$ and $F_{\bi{q}} = 0$\\
1& otherwise,\\}
\end{equation}
where $C$ is the complex conjugation operator. The second case above, although never strictly
encountered with realistic objects, expresses the fact that the rank of the projection effectively
decreases by one whenever $|\rho_{\bi{q}}| \gg F_{\bi{q}}$ (projection onto the origin, rather
than a circle). 
This case is included because small (unmeasured) moduli at large spatial frequencies are ubiquitous in
practice. Assuming for simplicity, however, that $Q_{\rm data}$ includes all $N$ Fourier samples, and that all the
moduli $F_{\bi{q}}$ are nonzero, the real-space form (\ref{realSpaceFourierProj}) 
of linearized Fourier modulus
projection is given by
\begin{equation}\label{linearFourierProj}
(\pi_{\rm F})_{\bi{r\, r^\prime}} = \frac{1}{2}\delta_{\bi{r\, r^\prime}}-
\frac{1}{2 N}\sum_{\bi{q}}
\exp{\rmi \left[2\phi_{\bi{q}}-\bi{q}\cdot (\bi{r+r^\prime)}\right]}\; C\;.
\end{equation}

To facilitate the comparison with random matrix theory, we express $\pi_{\rm F}$
in somewhat more compact notation. In the case of real-valued objects, where the operator $C$ can be dropped, we have
\begin{equation}\label{FourierProjReal}
\pi_{\rm F}=\frac{1}{2}(1+U^\mathrm{T} U)\;,
\end{equation}
where the unitary matrix $U$ is defined by
\begin{equation}
(U)_{\bi{q\, r}}=\frac{1}{\sqrt{N}}\exp{\rmi (\phi_{\bi{q}}-\bi{q}\cdot\bi{r})}\;.
\end{equation}
Since for real objects $\phi_{\bi{-q}}=-\phi_{\bi{q}}$, the reality of $U^\mathrm{T} U$ follows from
\begin{equation}
(U)_{\bi{q\, r}}^*=(U)_{\bi{-q\, r}}\;.
\end{equation}
When representing complex-valued objects in a $2N$-dimensional real vector space,
a general operator $X+Y C$, where $X$ and $Y$ are ordinary matrices, has the
block structure,
\begin{equation}\label{blockForm}
\left(
\begin{array}{cc}
{\rm Re}(X+Y) & {\rm Im}(Y-X)\\
{\rm Im}(X+Y) & {\rm Re}(X-Y)
\end{array}
\right)\;,
\end{equation}
and the projection takes the form
\begin{equation} \label{FourierProjComplex}
\pi_{\rm F}=\frac{1}{2}\left(
\begin{array}{cc}
1+{\rm Re}(U^\mathrm{T} U) & {\rm Im}(U^\mathrm{T} U)\\
{\rm Im}(U^\mathrm{T} U) & 1-{\rm Re}(U^\mathrm{T} U)
\end{array}
\right)\;.
\end{equation}

The linearized difference map is defined by,
\begin{equation}
d(\eta) :=\lim_{\epsilon\to 0}\frac{D(\rho^{\rm sol}+\epsilon\eta) - \rho^{\rm sol}}{\epsilon}\;,
\end{equation}
and is obtained by replacing the projections $\Pi_1$ and $\Pi_2$ in $D$ by their linearized counterparts,
$\pi_1$ and $\pi_2$.
The fixed points of $d$ form a linear space that
can be directly constructed from the kernels of $\pi_1$ and $\pi_2$ as follows. From (\ref{fixPtProj}) and
$\rho^* = \rho^{\rm sol} + \eta^*$ we have
\numparts
\begin{eqnarray}
\pi_1\circ\left[(1+\gamma_2)\pi_2 - \gamma_2\right](\eta^*) = & 0\label{kernel1}\\
\pi_2\circ\left[(1+\gamma_1)\pi_1 - \gamma_1\right](\eta^*) = & 0\label{kernel2}\;.
\end{eqnarray}
\endnumparts
By forming suitable linear combinations of (\ref{kernel1}) and (\ref{kernel2}), and making use
of the idempotency of $\pi_1$ and $\pi_2$, one can show
\begin{eqnarray}
\pi_1(\eta^*) = & 0\\
\pi_2(\eta^*) = & 0\;.
\end{eqnarray}
Thus,
\begin{equation}\label{diffKernel}
\ker{(d-1)} = \ker{\pi_1}\cap\ker{\pi_2}:=\ker{\pi_\perp}\;,
\end{equation}
where $\pi_\perp$ is the projection onto the orthogonal complement of the linear space
of fixed points. Since it is the projection of the iterates $d^n(\eta)$ by $\pi_\perp$ that should converge
to zero, we consider in the next Section the operator
\begin{equation}\label{dPerp}
d_\perp :=\pi_\perp\circ d \circ\pi_\perp\;.
\end{equation}
Applying DeMorgan's law (for Euclidean subspaces) to (\ref{diffKernel}), we obtain
\begin{equation}\label{DeMorgan}
\ker{(1-\pi_1)}\cup\ker{(1-\pi_2)}=\ker{(1-\pi_\perp)}\;.
\end{equation}
By using (\ref{diffKernel}) and (\ref{DeMorgan}) one easily verifies a set of identities needed in the
next Section:
\begin{equation}\label{perpProjIdentity}
\pi_i\circ\pi_\perp = \pi_\perp\circ\pi_i = \pi_i\qquad(i = 1, 2)\;.
\end{equation}
Finally, since we confine ourselves to the overconstrained case where
\begin{equation}
\ker{(1-\pi_1)}\cap\ker{(1-\pi_2)}=\{0\}\;,
\end{equation}
it follows from (\ref{DeMorgan}) that
\begin{equation}
\dim{\ker{(1-\pi_1)}}+\dim{\ker{(1-\pi_2)}}=\dim{\ker{(1-\pi_\perp)}}\;,
\end{equation}
or equivalently,
\begin{equation}\label{perpProjTrace}
\Tr{\pi_1} + \Tr{\pi_2}=\Tr{\pi_\perp}\;.
\end{equation}

\section{Convergence of the difference map}
With linearization, convergence to the space of fixed points corresponds to the vanishing of the
norm of the iterates $(d_\perp)^n(\eta)$, where the map $d_\perp$ was defined in (\ref{dPerp}).
This in turn is guaranteed by a suitable bound on an induced matrix norm of $d_\perp$. With the
standard Euclidean norm on the vector space, the natural choice for the induced matrix 
norm is the spectral norm:
\begin{equation}
\|d_\perp\|_\infty:=\max_{|\eta|=1} |d_\perp(\eta)|\;.
\end{equation}
There are two reasons, however, why the spectral norm may not be the appropriate choice. First,
it presents a serious challenge to calculate, even in the limit $N\to\infty$. Second, although being
able to achieve $\|d_\perp\|_\infty<1$ by suitable choice of parameters guarantees convergence,
it may not be true that minimizing this bound also maximizes the average rate of convergence.
An alternative is suggested by the interpretation of $(\|d_\perp\|_\infty)^2$ as the standard
$\infty$-norm applied to the eigenvalues of the matrix ${d_\perp}^\dagger\circ d_\perp$. 
In particular, the 1-norm
of the (non-negative) eigenvalues has the advantage that it is easily calculated:
\begin{equation}
(\|d_\perp\|_1)^2 := \frac{1}{N}\Tr{{d_\perp}^\dagger\circ d_\perp}\;.
\end{equation}
The norm $\|d_\perp\|_1$, called the Frobenius norm, 
gives the root-mean-square change in the Euclidean length of uniformly
sampled unit vectors when acted on by $d_\perp$. 
Although a bound on $\|d_\perp\|_1$ does not guarantee convergence, its minimization makes more
sense in that it reflects a property of all the eigenvalues rather than just the top of the
spectrum. This feature, and its tractable form, leads us to adopt the Frobenius norm as the vehicle
for optimizing the difference map; the subscript 1 will subsequently be dropped.

A straightforward calculation, making repeated use of (\ref{perpProjIdentity}), gives
\begin{eqnarray}\label{FrobeniusNorm}\fl
\|d_\perp\|^2 &= t_\perp +
2\beta\left[\gamma_1 t_2 - \gamma_2 t_1 + (\gamma_2-\gamma_1)t_{12}\right]+\nonumber \\
\fl &\beta^2\left[\gamma_1^2 t_2 + \gamma_2^2 t_1 +
\left(2+2(\gamma_1+\gamma_2)-(\gamma_1-\gamma_2)^2\right)t_{12}-
2(1+\gamma_1)(1+\gamma_2)t_{1212}\right]
\; ,
\end{eqnarray}
where
\numparts
\begin{eqnarray}
t_1 &:= (\Tr{\pi_1})/N\\
t_2 &:= (\Tr{\pi_2})/N\\
t_\perp &:= (\Tr{\pi_\perp})/N = t_1+t_2\label{tracePerpProj}\\
t_{12} &:= (\Tr{\pi_1\circ\pi_2})/N\\
t_{1212} &:= (\Tr{\pi_1\circ\pi_2\circ\pi_1\circ\pi_2})/N\;.
\end{eqnarray}
\endnumparts
Since (\ref{FrobeniusNorm}) is quadratic in the $\gamma_i$ and non-negative, it has a unique
minimum where the partial derivatives $\partial(\|d_\perp\|^2)/\partial\gamma_i$ vanish. The optimal
$\gamma_i$ are given by the following expressions:
\numparts
\begin{eqnarray}
\gamma_1^{\rm opt} &=-\frac{
(t_1-t_{12})(t_2-t_{1212})-(t_{12}-t_{1212})(t_1-2 t_{12}+t_{1212})\beta}
{\left[ t_1(t_2-t_{12})+t_{12}(2 t_{1212}-t_2)-t_{1212}^2\right]\beta}\label{gamma1Opt}\\
\bs
\gamma_2^{\rm opt} &=\frac{
(t_2-t_{12})(t_1-t_{1212})-(t_{12}-t_{1212})(t_2-2 t_{12}+t_{1212})\beta}
{\left[ t_2(t_1-t_{12})+t_{12}(2 t_{1212}-t_1)-t_{1212}^2\right]\beta}\label{gamma2Opt}\;.
\end{eqnarray}
\endnumparts
The Frobenius norm (\ref{FrobeniusNorm}) evaluated for these parameter values is
\begin{eqnarray}\label{evalFrobeniusNorm}
\fl
\|d_\perp\|^2 = t_\perp\nonumber\\
\fl+\left\{
(t_1-t_{12})(t_2-t_{12})(2 t_{1212}-t_1-t_2)+
2(t_2-t_1)(t_{12}-t_{1212})^2\beta\right.\nonumber\\
\left.\fl \quad+\left[(t_{12}-t_{1212})(2 t_1 t_2-3 t_{12}(t_1+t_2)+4 t_{12}^2+t_{1212}(t_1+t_2-2 t_{12})
\right]\beta^2 \right\}/\nonumber\\
\left\{t_1 t_2-t_{12}(t_1+t_2-2 t_{1212})-t_{1212}^2\right\}\;.
\end{eqnarray}
When the linearized constraint subspaces are orthogonal, as was assumed in reference\cite{elser},
the product traces $t_{12}$ and $t_{1212}$ vanish and one recovers the result
$\gamma_1^{\rm opt} = -1/\beta$, $\gamma_2^{\rm opt} = 1/\beta$. Using
(\ref{tracePerpProj}) in (\ref{evalFrobeniusNorm}), one obtains $\|d_\perp\|^2 = 0$ for this case.
This simple result could have been anticipated since for the given values of the $\gamma_i$ and
orthogonal subspaces ($\pi_1\circ\pi_2 = 0$) one has the simplification
$d=1-\pi_1-\pi_2=1-\pi_\perp$, and hence $d_\perp=0$.

\section{Average traces of projection products}\label{aveTrace}
Because our convergence measure $\|d_\perp\|^2$ is linear in the traces $t_{12}$ and $t_{1212}$,
an average of these traces with respect to a particular class of objects will provide the
optimum average-case convergence. The projections that appear in these traces 
can be expressed as $\pi=O^\mathrm{T}PO$, where $P$ is a fixed 
canonical projection and $O$ is an orthogonal matrix. The actual
matrices $O$ that apply to the projections considered earlier comprise subensembles in the
space of orthogonal matrices. For $\pi_{\rm S}$ (object support) $O$ is a
permutation matrix, while for $\pi_{\rm F}$ (Fourier modulus) $O$ belongs to a continuous
family of orthogonal matrices defined implicitly by (\ref{FourierProjReal}) or (\ref{FourierProjComplex}).
Traces of products of $\pi_1={O_1}^\mathrm{T}P_1O_1$ and $\pi_2={O_2}^\mathrm{T}P_2O_2$
only depend on the relative orthogonal matrix $O=O_1{O_2}^\mathrm{T}$. 
In the absence of detailed knowledge of the appropriate subensemble for $O$, we might
consider the average over the complete, group invariant ensemble:
\numparts
\begin{eqnarray}
\langle t_{12}\rangle &:=\frac{1}{N}\langle\Tr{P_1 O^\mathrm{T} P_2 O}\rangle\\
\langle t_{1212}\rangle &:=\frac{1}{N}\langle\Tr{P_1 O^\mathrm{T} P_2 O P_1 O^\mathrm{T} P_2 O}\rangle\;.
\end{eqnarray}
\endnumparts
There are highly sophisticated methods \cite{creutz, samuel} for evaluating 
these standard random matrix averages
which rely solely on symmetry.
Here we quote the limiting forms for large $N$; details of the calculation are given in the Appendix:
\numparts
\begin{eqnarray}
\lim_{N\to\infty}\langle t_{12}\rangle &= t_1 t_2\label{t12}\\
\lim_{N\to\infty}\langle t_{1212}\rangle &= t_1^2 t_2 + t_1 t_2^2 - t_1^2 t_2^2\label{t1212}\;.
\end{eqnarray}
\endnumparts
In the remainder of this Section these formulas are
compared with explicit averages over actual {\em object ensembles} (and
$\pi_1 = \pi_{\rm S}$, $\pi_2 = \pi_{\rm F}$).
If the random matrix
theory results (\ref{t12}, \ref{t1212}) are reproduced,
then the problem of choosing optimal parameters is greatly simplified since (\ref{FrobeniusNorm})
will then only depend
on the dimensionalities of the constraint subspaces ($t_1$ and $t_2$).

We begin by calculating the traces of the two projections. From (\ref{suppProj}) and
(\ref{linearSuppProj}) we have
\begin{equation}
t_{\rm S} = (\Tr{\pi_{\rm S}})/N = |S|/N := \sigma\;,
\end{equation}
where $\sigma$ is the fraction of pixels or voxels in the support $S$. 

For traces involving the linearized Fourier
modulus projection (\ref{linearFourierProj}) we need to understand the action of the conjugation
operator $C$. In the case of real-valued objects all traces may be evaluated in a real vector space of
dimension $N$, and $C$ can be replaced by the identity. For complex-valued objects traces are
evaluated in a $2N$-dimensional real vector space where the representation of a general operator
$X+Y C$ has the block form (\ref{blockForm}).
Normalizing the matrix trace of (\ref{blockForm}) by an extra factor of $\frac{1}{2}$ gives us the
identity
\begin{equation}\label{traceIdentity}
\Tr{(X + Y C)} := \Tr{({\rm Re}X)}\;,
\end{equation}
where the trace on the right is an ordinary matrix trace. Another useful identity
involving $C$ is
\begin{equation}\label{prodIdentity}
(X+Y C)(W+Z C) = (X W + Y Z^*) + (X Z + Y W^*)C\;,
\end{equation}
where $X$, $Y$, $W$, and $Z$ are ordinary matrices. 

Use of identity (\ref{traceIdentity}) on (\ref{linearFourierProj}) gives us immediately
 \begin{equation}\label{traceFcomplex}
t_{\rm F} = (\Tr{\pi_{\rm F}})/N = \frac{1}{2}\;
\end{equation}
in the case of complex-valued objects.
For real-valued objects we have
\begin{equation}\label{traceFreal}
(\Tr{\pi_{\rm F}})/N = \frac{1}{2} -
\frac{1}{2 N^2}\sum_{\bi{q}}
\exp{(\rmi 2\phi_{\bi{q}})}\left[{\textstyle \sum_{\bi{r}}}\exp{(-\rmi 2\bi{q}\cdot \bi{r})}\right]\;.
\end{equation}
The inner sum in (\ref{traceFreal}) vanishes unless the components of $\bi{q}$ satisfy
$q_i \equiv 0\pmod{\pi m_i}$, and there are at most 4 (two dimensions) or 8 (three dimensions)
such points. In more concrete terms, the second term in (\ref{traceFreal}) subtracts the number of samples
$\bi{q}$ which do not possess a continuous phase when the object is real. This correction is negligible
in the $N\to\infty$ limit and we recover the complex-object result (\ref{traceFcomplex}).

We consider next the first product trace computed earlier using random matrix theory. In the
complex-object case, identities (\ref{traceIdentity}) and (\ref{prodIdentity}) give us immediately
\begin{equation}
t_{\rm SF} = (\Tr{\pi_{\rm S}\circ\pi_{\rm F}})/N = \frac{\sigma}{2}\;,
\end{equation}
in agreement with the random matrix result (\ref{t12}), even without averaging. 
For real-valued objects we have
\begin{equation}\label{traceSFreal}
(\Tr{\pi_{\rm S}\circ\pi_{\rm F}})/N = \frac{\sigma}{2} -
\frac{1}{2 N}\sum_{\bi{q}}
\exp{(\rmi 2\phi_{\bi{q}})}\,\Sigma_{-2\bi{q}}\;,
\end{equation}
where
\begin{equation}
\Sigma_{\bi{q}} := \frac{1}{N}\sum_{\bi{s}\in S}\exp{(\rmi \bi{q}\cdot \bi{s})}
\end{equation}
defines the Fourier transform of the support's characteristic function. 

The second term in (\ref{traceSFreal}) is a correction to the naive random matrix theory result and in general
does not vanish. To check whether it vanishes in an average-case sense, we consider two object
ensembles distinguished by the nature of their support. If the support is a known fixed set $S$, say
a rectangle or disk in two dimensions, then $\Sigma_{-2\bi{q}}$ is a constant in the ensemble average
of (\ref{traceSFreal}) and we can focus on the average of $\exp{(\rmi 2\phi_{\bi{q}})}$.
The minimal {\it a priori} knowledge of the object in real-space is that the individual pixel (or voxel)
values $\rho_{\bi{r}}$ in the support are independently distributed according to a known distribution. For this
object ensemble it is straightforward to calculate the distribution of the Fourier transform
components $\rho_{\bi{q}}$. Excepting a fraction of samples $\bi{q}$ which vanishes as $N\to \infty$
(near $\bi{q} = 0$ or associated with facets in the support), the complex numbers $\rho_{\bi{q}}$
are isotropic Gaussian random variables with mean zero and corrections that decay in the limit $N\to \infty$.
Consequently, the average of the phases $\exp{(\rmi 2\phi_{\bi{q}})}$ in (\ref{traceSFreal}) vanishes
in the same limit, as does the correction (since $|\Sigma_{\bi{q}}|\leq\sigma$, independent of $N$).

In the case of an atomic support, although the actual support $S$ is not known, one does know that
$S$ is the union of a given number $M$ of very compact sets associated with atoms. In the equal-atom
ensemble there is an approximate relationship
\begin{equation}\label{formFactorRelation}
\Sigma_{\bi{q}}\approx\frac{A(\bi{q})}{\sqrt{N}}\rho_{\bi{q}}\;,
\end{equation}
where $A(\bi{q})$ is a smooth real-valued function (form factor) independent of $N$. The approximation in
(\ref{formFactorRelation}) arises from shifts (by fractions of voxels) of the atomic centers relative to their individual supports. To bound the sum in (\ref{traceSFreal}) we replace each term by its modulus and
use (\ref{formFactorRelation}):
\begin{equation}\label{bound}
\frac{1}{2 N}\left|\sum_{\bi{q}}
\exp{(\rmi 2\phi_{\bi{q}})}\,\Sigma_{-2\bi{q}}\right|  <
\frac{1}{2 N}\sum_{\bi{q}} \frac{|A(2\bi{q})|}{\sqrt{N}}|\rho_{2\bi{q}}|\;.
\end{equation}
The root-mean-square of the modulus $|\rho_{2\bi{q}}|$ can be calculated in the 
standard atomic ensemble where $M$ atom centers are uniformly and independently distributed; the result,
$|\rho_{2\bi{q}}|_{\rm rms}=\Or{(\sqrt{M/N})}$, ($2\bi{q}\neq 0$) shows that the correction is negligible
for $N\to \infty$ and fixed resolution ($M/N = {\rm const}$).
We conclude that explicit averaging in the standard equal-atom ensemble,
as in the fixed support ensemble,
reproduces the random matrix theory result (\ref{t12}):
\begin{equation}
\langle t_{\rm SF} \rangle = \frac{\sigma}{2}\;.
\end{equation}

Finally, we turn to the trace of the product of four projection operators. For complex-valued objects
we obtain
\begin{equation}\label{traceSFSFcomplex}
\fl
t_{\rm SFSF}=(\Tr{\pi_{\rm S}\circ\pi_{\rm F}\circ\pi_{\rm S}\circ\pi_{\rm F}})/N = \frac{\sigma}{4} +
\frac{1}{4 N}\sum_{\bi{q},\bi{q^\prime}}
\exp{\rmi 2(\phi_{\bi{q}}-\phi_{\bi{q^\prime}})}\,\Sigma_{\bi{q^\prime}-\bi{q}}^2\;,
\end{equation}
and for real objects:
\begin{eqnarray}\label{traceSFSFreal}
\fl
t_{\rm SFSF}=(\Tr{\pi_{\rm S}\circ\pi_{\rm F}\circ\pi_{\rm S}\circ\pi_{\rm F}})/N = \frac{\sigma}{4} -
\frac{1}{2 N}\sum_{\bi{q}}
\exp{(\rmi 2\phi_{\bi{q}})}\,\Sigma_{-2\bi{q}}\nonumber\\
+\frac{1}{4 N}\sum_{\bi{q},\bi{q^\prime}}
\exp{\rmi 2(\phi_{\bi{q}}+\phi_{\bi{q^\prime}})}\,\Sigma_{-\bi{q^\prime}-\bi{q}}^2\;.
\end{eqnarray}
We recognize the single sum in (\ref{traceSFSFreal}) as the sum in (\ref{traceSFreal}), which was
already argued to be negligible when averaged over the standard object ensembles with fixed or atomic support. 
However, the averages of the double sums for the fixed support ensemble, 
both in (\ref{traceSFSFcomplex}) and (\ref{traceSFSFreal}),
can deviate from the random matrix theory result (\ref{t1212}). In fact, the latter
is recovered in both cases by retaining only the diagonal terms in the double sums
($\bi{q}=\bi{q}^\prime$ in (\ref{traceSFSFcomplex}), 
$\bi{q}=-\bi{q}^\prime$ in (\ref{traceSFSFreal})). Since the contributions of the near-diagonal terms
is not diminished in the $N\to\infty$ limit, there will be non-negligible corrections that depend on details
of the ensemble (support shape, distribution of object values, etc.). The most we can conclude for the
{\em fixed support ensemble} without these details is,
\begin{equation}
\langle t_{\rm SFSF}\rangle = \frac{\sigma}{4}+\Or{(\sigma^2)}\;,
\end{equation}
a result that may be useful for small object support.

The neglect of non-diagonal terms in the double sum may be easier to justify rigorously for the
atomic support ensemble. Considering only the case of real atoms (\ref{traceSFSFreal}), we first
use (\ref{formFactorRelation}) to write
\begin{equation}
\Sigma_{-\bi{q^\prime}-\bi{q}}^2\approx
\frac{|A(-\bi{q^\prime}-\bi{q})|^2}{N}|\rho_{-\bi{q^\prime}-\bi{q}}|^2
\exp{(-2\rmi \phi_{\bi{q^\prime}+\bi{q}})}\;.
\end{equation}
Now, since $|\rho_{-\bi{q^\prime}-\bi{q}}|^2_{\rm rms} = \Or{(M/N)}$ for $\bi{q}+\bi{q}^\prime\neq 0$,
the non-diagonal terms have moduli which scale as $N^{-2}$ (at constant resolution). To argue that the
$\Or{(N^2)}$ non-diagonal terms in the double
sum have a negligible ensemble average, we consider the phase of a general
term:
\begin{equation}\label{triplet}
\exp{\rmi 2(\phi_{\bi{q}}+\phi_{\bi{q^\prime}}-\phi_{\bi{q^\prime}+\bi{q}})}\;.
\end{equation}
Atomic ensemble averages of {\em triplet phase invariants} such as (\ref{triplet})
are studied at length in the statistical theory of structure factors in crystallography\cite{giacovazzo}; 
they are found
to vanish in the limit of many atoms, or equivalently, $N\to\infty$ when the resolution is fixed. 
There are three exceptions
to this general result, when $\bi{q}=0$, $\bi{q^\prime}=0$, or $\bi{q^\prime}+\bi{q}=0$. The
latter corresponds to the diagonal contribution in the double sum of (\ref{traceSFSFreal}) which
we will retain. The two other cases lead to single sums whose moduli can be shown to be negligible
by arguments identical to those used for the sum in (\ref{bound}). In summary, this non-rigorous
argument leads us to believe that for the {\em atomic support ensemble}
\begin{equation}
\langle t_{\rm SFSF}\rangle = \frac{\sigma}{4}+\frac{\sigma^2}{4}\;,
\end{equation}
in agreement with the random matrix theory result (\ref{t1212}).

\section{Optimal parameter values}
The results of the previous Section show that random matrix averaging, of
the traces of products of projections, reproduces for the most part the averages over the actual
object ensembles encountered in standard phase retrieval applications (fixed and atomic support).
This greatly simplifies the expressions for the optimal $\gamma$ parameters
of the difference map (\ref{gamma1Opt}, \ref{gamma2Opt}), since all the traces are explicit functions of just
the traces $t_{\rm S}=\sigma$ and $t_{\rm F}=\frac{1}{2}$:
\numparts
\begin{eqnarray}
\gamma_{\rm S}^{\rm opt}&=-\frac{4+(2+\beta)\sigma+\beta\sigma^2}{\beta(4-\sigma+\sigma^2)}\\
\bs
\gamma_{\rm F}^{\rm opt}&=\frac{6-2\sigma-\beta(2-3\sigma+\sigma^2)}{\beta(4-\sigma+\sigma^2)}\;.
\end{eqnarray}
\endnumparts
Following the discussion of uniqueness in Section 2, phase retrieval with a 
support constraint is well posed when $t_{\rm S}+t_{\rm F}<1$, which requires that the support fraction
satisfies $0<\sigma<\frac{1}{2}$.
As noted in Section 5, with the fixed support ensemble the average of $t_{\rm SFSF}$ (\ref{traceSFSFcomplex},
\ref{traceSFSFreal}) may differ from the random matrix average by terms proportional to $\sigma^2$,
and the coefficients of $\sigma^2$ in the above formulae will change accordingly. For the atomic support
ensemble the formulae as written are more reliable. Moreover, in the principal atomic support application,
crystallography, there is a large dimensionality asymmetry of the two constraint subspaces as expressed
by the statement $\sigma \ll \frac{1}{2}$. For this important application one may therefore take the
$\sigma\to 0$ limit:
\numparts
\begin{eqnarray}
\gamma_{\rm S}^{\rm opt}&\stackrel{\sigma\to 0}{=}
-\frac{1}{\beta}\\
\bs
\gamma_{\rm F}^{\rm opt}&\stackrel{\sigma\to 0}{=}
\frac{3-\beta}{2\beta}\label{gammaFlimit}\;.
\end{eqnarray}
\endnumparts
The Frobenius norm (\ref{evalFrobeniusNorm}) in this limit is given by
\begin{equation}\label{normlimit}
\|d_\perp\|^2 \stackrel{\sigma\to 0}{=}
\frac{\sigma}{8}(3+2\beta+3\beta^2)\;.
\end{equation}

The numerical constants in (\ref{gammaFlimit}) and (\ref{normlimit}) can be traced to the numerical
value of $t_{\rm F}=\frac{1}{2}$, the only trace with a finite value in the $\sigma\to 0$ limit.
Recalling the discussion in Section 3.3 and the interpretation of $t_{\rm F}N$ as the dimensionality
of the Fourier modulus constraint subspace, the exact value $\frac{1}{2}$ is certainly an oversimplification.
As a first improvement in the evaluation of $\Tr{\pi_{\rm F}}$ one should exclude contributions from
samples $\bi{q}$ for which the corresponding modulus $F_{\bi{q}}$ is negligible. In crystallography,
for example, only the samples within an ellipsoid about the origin have non-negligible $F_{\bi{q}}$, effectively
reducing the rank of the Fourier modulus projection by a fraction given approximately (in three dimensions) by
$\pi/6\approx 0.52$; thus $t_{\rm F}\approx 0.26$. Given the uncertainty in $t_{\rm F}$ we
give formulas for the optimal $\gamma$ parameters for general $t_{\rm F}$:
\numparts
\begin{eqnarray}
\gamma_{\rm S}^{\rm opt}&\stackrel{\sigma\to 0}{=}
-\frac{1}{\beta}\\
\bs
\gamma_{\rm F}^{\rm opt}&\stackrel{\sigma\to 0}{=}
\frac{1+t_{\rm F}(1-\beta)}{\beta}\;.
\end{eqnarray}
\endnumparts

\ack
The author thanks Orlando Alvarez and Piet Brouwer for help with the random matrix averages,
and the referee for exceptionally constructive criticism.
Additional thanks go to the Center for Experimental and Constructive Mathematics, where the
calculations were performed, and the Aspen Center for Physics, where this article was prepared.
This work was supported by the National Science Foundation under grant ITR-0081775.

\appendix
\section*{Appendix}
\setcounter{section}{1}
Averages over products of matrix elements of $N\times N$ orthogonal matrices $O$
($N>1$), such as
\begin{equation}\label{twoO}
X^{i\,j}_{p\,q}=\langle O^i_{\;p} O^j_{\;q}\rangle\;,
\end{equation}
can always be expressed as tensor products of Kronecker-$\delta$'s having the correct
invariance properties. The latter follow from the invariance of the group integration measure
with respect to the change of variables $O\to L O R$, where $L$ and $R$ are independent
orthogonal matrices. For the example above, the statement of invariance becomes
\begin{equation}
L^i_{i^\prime}L^j_{j^\prime}\;X^{i^\prime\, j^\prime}_{p^\prime\, q^\prime}\;
R^{p^\prime}_{p}R^{q^\prime}_{q}=X^{i\,j}_{p\,q}\;,
\end{equation}
with general solution
\begin{equation}\label{sol2}
X^{i\,j}_{p\,q}=a\,\delta^{ij}\delta_{pq}\;,
\end{equation}
where $a$ is an undetermined constant. By contracting the indices $i$ and $j$, whereupon
the left hand side of (\ref{sol2}) 
corresponds to $\langle(O^\mathrm{T}O)_{pq}\rangle=\delta_{pq}$, consistency requires
\begin{equation}
a=1/N\;.
\end{equation}

Using this method for the product of four matrix elements,
\begin{equation}\label{fourO}
X^{i\,j\,k\,l}_{p\,q\,r\,s}=\langle O^i_{\;p} O^j_{\;q} O^k_{\;r} O^l_{\;s}\rangle\;,
\end{equation}
we find
\begin{eqnarray}\label{sol4}
X^{i\,j\,k\,l}_{p\,q\,r\,s}=&a_1&
(\delta^{ij}\delta^{kl}\delta_{pq}\delta_{rs}+
\delta^{ik}\delta^{jl}\delta_{pr}\delta_{qs}+
\delta^{il}\delta^{jk}\delta_{ps}\delta_{qr})+\\ \nonumber
&a_2&
(\delta^{ij}\delta^{kl}\delta_{pr}\delta_{qs}+
\delta^{ij}\delta^{kl}\delta_{ps}\delta_{qr}+
\delta^{ik}\delta^{jl}\delta_{pq}\delta_{rs}+\\ \nonumber
& &
\delta^{ik}\delta^{jl}\delta_{ps}\delta_{qr}+
\delta^{il}\delta^{jk}\delta_{pq}\delta_{rs}+
\delta^{il}\delta^{jk}\delta_{pr}\delta_{qs})\;.
\end{eqnarray}
The reduction to two undetermined constants made use
of permutation symmetries. Contracting indices $k$ and $l$ in (\ref{fourO})
gives our previous result for (\ref{twoO}) multiplied by $\delta_{rs}$, which can then
be compared with the result of contracting the same pair of indices in
(\ref{sol4}). Consistency yields two equation with the solutions
\begin{eqnarray}
a_1&=&\frac{N+1}{(N-1)N(N+2)}\\
a_2&=&\frac{-1}{(N-1)N(N+2)}\;.
\end{eqnarray}

The random matrix averages of the traces of products of projections are obtained
by contracting the tensor expressions (\ref{sol2}) and (\ref{sol4}) with the projections
$P_1$ and $P_2$:
\begin{eqnarray}\fl
\langle\Tr{P_1 O^\mathrm{T} P_2 O}\rangle&=&a(\Tr{P_1})(\Tr{P_2})\\ \fl\label{sol4P}
\langle\Tr{P_1 O^\mathrm{T} P_2 O P_1 O^\mathrm{T} P_2 O}\rangle&=&
a_1(\Tr{P_1})(\Tr{P_2})(1+\Tr{P_1}+\Tr{P_2})+\\ \fl\nonumber
& &a_2(\Tr{P_1})(\Tr{P_2})(3+\Tr{P_1}+\Tr{P_2}+(\Tr{P_1})(\Tr{P_2}))\;.
\end{eqnarray}
After taking the limit $N\to\infty$ in (\ref{sol4P}) we obtain the results quoted in Section \ref{aveTrace}.

\end{document}